\documentclass[final]{amsart}
\usepackage{amssymb,amscd}

\numberwithin{equation}{section}
\newtheorem {Theorem} {Theorem}
\newtheorem{Proposition} [Theorem] {Proposition}
\newtheorem{Lemma} [Theorem] {Lemma}
\newtheorem {Corollary} [Theorem] {Corollary}
\theoremstyle{remark}
\newtheorem* {Remark} {Remark}

\def \re        {relative equilibrium}
\def \rea       {relative equilibria}

\def \inv       {^{-1}}

\def \prc       {\Phi_{\C}}
\def \prr       {\Phi_{\R}}
\def \eps       {\varepsilon}
\def\sympquot#1#2{/\!\!/\!_#1 #2}
\def\restr#1{\vrule height1.2ex width.3pt
               depth1.4ex\lower1.0ex\hbox{\scriptsize $\,#1$}} 
\def\smallrestr#1{\vrule height0.8ex width.2pt
               depth0.8ex\lower0.6ex\hbox{\scriptsize $\,#1$}}

\newcommand{\fg} {{\mathfrak g}}
\newcommand{\fh} {{\mathfrak h}}
\newcommand{\fk} {{\mathfrak k}}
\newcommand{\fl} {{\mathfrak l}}
\newcommand{\fm} {{\mathfrak m}}
\newcommand{\ft} {{\mathfrak t}}

\newcommand{\so} {\mathfrak {so}}

\def    \cQ     {{\mathcal Q}}
\def    \cO     {{\mathcal O}}

\def    \C      {{\mathbb C}}

\def    \R      {{\mathbb R}}

\def\,{\mskip\thinmuskip}
\def\Coad{\textrm{Coad}}
\def\Ad{\textrm{Ad}}

\begin{document}

\title{Openness of momentum maps and persistence of extremal relative 
  equilibria}

\author{James Montaldi}
\author{Tadashi Tokieda}

\address{Department of Mathematics, UMIST, PO Box 88, Manchester M60 1QD, UK}
\email{j.montaldi@umist.ac.uk}

\address{D\'epartement de Math\'ematiques, Universit\'e de 
Montr\'eal, C.P.~6128, 
succ.~Centre-Ville, Montr\'eal H3C 3J7, Canada}
\email{tokieda@dms.umontreal.ca} 

\date{January, 2002}

\thanks{The research of JM was partially supported by the European
  Union through the Research Training Network MASIE}

\begin{abstract}
  We prove that for every proper Hamiltonian action of a Lie group $G$
  in finite dimensions the momentum map is locally $G$-open relative
  to its image (i.e.\ images of $G$-invariant open sets are open).  As
  an application we deduce that in a Hamiltonian system with
  continuous Hamiltonian symmetries, extremal relative equilibria
  persist for every perturbation of the value of the momentum
  map, provided the isotropy subgroup of this value is compact.  We
  also demonstrate how this persistence result applies to an example
  of ellipsoidal figures of rotating fluid, and provide an example
  with plane point vortices which shows how the compactness assumption
  is related to persistence.

  \medskip

\end{abstract}

\maketitle


\section{Introduction}

In a Hamiltonian system, nondegenerate equilibria are isolated; in
particular, they do not \emph{persist} from one energy level to nearby
levels.  In this paper we prove that in a symmetric Hamiltonian
system, every extremal relative equilibrium persists to nearby levels
of the momentum map, provided the isotropy subgroup of its momentum
value is compact.  The crucial ingredient in the proof is a
generalisation of Sjamaar's result on the openness of momentum maps.

Let $M$ be a symplectic manifold with a proper and symplectic action
of a connected Lie group $G$ and let $h$ be a $G$-invariant
Hamiltonian.  We suppose that the action of $G$ is Hamiltonian, in
that it is infinitesimally generated by a momentum map $\Phi: M \to
\fg^{\ast}$.  A trajectory of the Hamiltonian vector field $X_h$ of
$h$ is a \emph{relative equilibrium} if its image in the orbit space
$M/G$ is a single point; such a \re\ with momentum value $\mu$ is
\emph{extremal} if its image in the reduced space
$\Phi^{-1}(\mu)/G_{\mu}$ is a local extremum for the reduced
Hamiltonian.  $G_{\mu}$ denotes the isotropy subgroup of $\mu$ for the
possibly \emph{modified coadjoint action} of $G$ on $\fg^{\ast}$
(Section~\ref{sec:modified}).

Let $f:X\to Y$ be a $G$-equivariant map.  We say that $f$ is {\it $G$-open\/}
if the image of any $G$-invariant open set is open in $Y$.  This is
equivalent to the orbit map $\bar f:X/G\to Y/G$ being open.

Suppose that $\gamma$ is an extremal \re\ with momentum value $\mu$.
We wish to prove that under certain hypotheses such a \re\ persistes
to \emph{all} nearby values of the momentum map.  However, since this
is really a local result in the phase space $M$, the question arises
as to what is meant by all nearby values.  If the momentum map $\Phi$
is proper then a result of Sjamaar \cite{Sj} (see also \cite{LMTW})
says that the momentum map is $G$-open, so that images of
$G$-invariant neighbourhoods of $\gamma$ are open in the image of
$\Phi$, and the phrase `all nearby values' means just that: a full
$G$-invariant neighbourhood of $\mu$ in $\Phi(M)$.  On the other hand,
if $\Phi$ is not proper, then it may not be $G$-open (for an example,
see \cite[Example 3.10]{KL97}), so the image of a $G$-invariant
neighbourhood of $\gamma$ may not be open in $\Phi(M)$.  Our first
result shows that there is always a $G$-invariant neighbourhood $U_0$
of $\gamma$ restricted to which the momentum map $U_0\to\Phi(U_0)$
\emph{is} $G$-open relative to its image.  This neighbourhood $U_0$ is
a tubular neighbourhood of the group orbit containing $\gamma$ whose
existence follows from the Marle-Guillemin-Sternberg normal form for
symplectic group actions (Section \ref{sec:reducetocompact}).  
Of course, if $\Phi$ is proper one can take $U_0=M$.

\begin{Theorem}\label{thm:openness}
  Let $M$ be a symplectic manifold with a proper Hamiltonian action of
  a connected Lie group $G$ and a momentum map $\Phi : M \to
  \fg^{\ast}$.  Suppose that $x\in M$ has momentum value $\mu =
  \Phi(x)$ whose isotropy subgroup $G_{\mu}$ for the modified
  coadjoint action is compact.  Then there exists a $G$-invariant
  neighbourhood $U_0$ of $x$ such that the restriction
  $\Phi\restr{U_0}:U_0\to\Phi(U_0)$ is $G$-open, where $\Phi(U_0)$ is
  given the subspace topology induced from $\fg^{\ast}$.
\end{Theorem}

This result then allows us to state the persistence theorem for
extremal \rea.

\begin{Theorem} \label{thm:persistence}
  Let $M,G,\Phi,x,\mu$ and $U_0$ be as in Theorem \ref{thm:openness},
  let $h\in C^{\infty}(M)$ be a $G$-invariant Hamiltonian, and suppose
  that $\gamma$ is an extremal relative equilibrium for the given
  Hamiltonian system, with $x\in\gamma$.  Then there exists a
  $G$-invariant neighbourhood $V$ of $\mu$ in $\Phi(U_0)$ such that
  for every $\mu'\in V$, there is a \re\ in $\Phi\inv(\mu')\cap U_0$.
\end{Theorem}

When the $G$-action is trivial, $\Phi\inv(\mu') = M$, so the theorem
becomes trivial, too.  Theorems in Hamiltonian systems often have
natural generalisations to those in symmetric Hamiltonian systems,
when a group action is thrown in.  Theorem~\ref{thm:persistence} is an
instance of a theorem in the latter that has no nontrivial
specialisation in the former.

In \cite{Mo97}, it was shown that extremal relative equilibria are
Lyapunov-stable relative to $G$.  Also in \cite{Mo97} appeared a
version of persistence, but the proof is incomplete.  The present
version is stronger, first because it does not require $G$ to be compact,
but just $G_\mu$, which actually suffices to reduce to the compact
case (Section~\ref{sec:reducetocompact}), second because it
proves persistence to a full neighbourhood of $\mu$.  To our
knowledge, Theorem~\ref{thm:persistence} is the first application of
this topological property of the $G$-openness of momentum maps to
problems of Hamiltonian dynamical systems.

After its proof in Section~\ref{sec:persistence},
Theorem~\ref{thm:persistence} is applied in Section~\ref{sec:ARB} to
the problem of ellipsoidal figures of rotating fluid (affine rigid
bodies).  In Section~\ref{sec:vortices} we check that the compactness
hypothesis on the isotropy subgroup $G_\mu$ is essential in
Theorem~\ref{thm:persistence} by analysing point vortices on the
plane.  Finally, in Section~\ref{sec:stages} we explain how reduction
by stages yields a partial persistence result even in the case of
noncompact momentum isotropy.  For complementary results on
persistence of relative equilibria for noncompact group actions, see
Wulff \cite{Wulff01}.


\section{Modification of coadjoint action}\label{sec:modified}

Theorem~\ref{thm:openness} does not assume the equivariance of the
momentum map $\Phi$ with respect to the standard coadjoint action.  However,
Souriau \cite{Souriau70} showed that the momentum map can always be
made equivariant by modifying the coadjoint action, as follows.

Let a Lie group $G$ act in a Hamiltonian manner on a connected
symplectic manifold $M$ with a momentum map $\Phi : M \to \fg^{\ast}$.
Define the cocycle $\theta:G\to\fg^{\ast}$ by
$$
\theta(g) = \Phi(g\cdot x)-\Coad_g(\Phi(x))
$$
(which is independent of the choice of $x\in M$); $\Coad_g=
\Ad_{g^{-1}}^{\ast}$ is our notation for the standard coadjoint action
of $g\in G$ on $\fg^{\ast}$.  The \emph{modified coadjoint action} is
$$
\Coad^\theta_g(\mu) = \Coad_g(\mu)+\theta(g), \qquad \mu\in\fg^*. 
$$
With respect to this shifted affine action $\Phi$ becomes
equivariant.  All the usual properties of standard coadjoint action
continue to hold for the modified actions \cite{Souriau70}: for
example the momentum map is Poisson for a suitably modified Poisson
structure on $\fg^{\ast}$, and the symplectic leaves of the modified
Poisson structure are the modified coadjoint orbits.

Throughout this paper the reduced space at $\mu\in \fg^{\ast}$ is
understood to be $\Phi^{-1}(\mu)/G_\mu$, where $G_\mu$ is the isotropy
subgroup of $\mu$ for the modified coadjoint action.


\section{Reduction to actions of compact groups}\label{sec:reducetocompact}

Theorem~\ref{thm:openness} does not assume the compactness of the
symmetry group $G$, but only the compactness of the isotropy subgroup
$G_{\mu}$.  The reduction to compact group actions is based on the
Marle-Guillemin-Sternberg normal form for symplectic actions and
momentum maps, which we now recall \cite{Marle, GS}.
%
%
Let again a connected Lie group $G$ act in a Hamiltonian manner on a
symplectic manifold $M$ with a momentum map $\Phi : M \to \fg^{\ast}$
and the corresponding cocycle $\theta: G \to \fg^{\ast}$
(Section~\ref{sec:modified}).  At $x\in M$, consider the four spaces
\begin{align*}
  T_0 &= T_x(G\cdot x)\cap \ker d\Phi(x) = T_x(G_{\mu}\cdot x), \\
  T_1 &= T_x(G\cdot x)/T_0, \\
  N_1 &= \ker d\Phi(x)/T_0, \\
  N_0 &= T_xM/(\, T_x(G\cdot x) + \ker d\Phi(x)).
\end{align*}
Since $\ker d\Phi_x$ is the symplectic complement to $T_x(G\cdot x)$, these 
spaces depend only on the $G$-action and not on the choice of $\Phi$.  Using 
the compactness of $G_x \subset G_{\mu}$, we can realise the quotients
$T_1$, $N_1$, $N_0$ as $G_x$-invariant subspaces of $T_xM$ satisfying
$$
T_0 \oplus T_1 = T_x(G\cdot x), \quad T_0 \oplus N_1 =
\ker d\Phi(x), \quad T_0 \oplus T_1 \oplus N_1 \oplus N_0 = T_xM
$$%
(the so-called Witt or Moncrief decomposition).  $N_1$ is the
\emph{symplectic slice\/} to the action at $x$. With respect to such a
decomposition, the symplectic form $\omega$ has the matrix
\begin{equation} \label{eq:matrix}
  [\omega]_x=\left[\begin{array}{cccc}
      0&0&0&A\cr 0&\omega_{T_1}&0&*\cr 0&0&\omega_{N_1}&*\cr -A^t&*&*&*
    \end{array}\right],
\end{equation}
where $\omega_{T_1}$ and $\omega_{N_1}$ are the restrictions of
$\omega$ to $T_1$ and $N_1$, $A$ is nondegenerate, and the $*$'s are
of no interest.  The Marle-Guillemin-Sternberg normal form theorem
states that in a $G$-invariant neighbourhood $U$ of $x$, the
symplectic $G$-action is isomorphic to that on
$G\times_{G_x}(\fm^*\times N_1)$, where $\fm^* = \fg_x^\circ \cap
\fg_\mu^*$, so that $\fg_\mu\simeq\fm\oplus\fg_x$ and
$\fg_\mu^*\simeq\fm^*\oplus\fg_x^*$.  ($\fg_x^\circ\cap\fg_\mu^*$ is
the annihilator of $\fg_x$ in $\fg_\mu^*$.)  The momentum map has the
explicit form
\begin{equation} \label{eq:MGS}
  \begin{array}{rcl}
    \Phi :G\times_{G_x} (\fm^*\times N_1) &\longrightarrow&  \fg^*\cr
    [g,\,\nu,\,v] &\longmapsto& \Coad_g^\theta(\mu+(\nu\oplus\Phi_{G_x}(v)). 
  \end{array}
\end{equation}

   Now we reduce the problem to the case where the whole group $G$ is
compact.  (This part of the argument is similar to the beginning of
Section 2 in \cite{LT}.)  Since the isotropy subgroup $G_{\mu}$ of $\mu$
is compact, we can choose a momentum map so that the restriction of
$\theta$ to $G_\mu$ vanishes (essentially by averaging \cite{Mo97}), so 
that for $g\in G_\mu$ we have $\Phi(g\cdot x) =
\Coad_g\Phi(x)$.  There is therefore an inner product on $\fg^*$,
invariant under $\Coad(G_{\mu})$, inducing a $G_{\mu}$-equivariant 
splitting $\fg = \fg_{\mu} \oplus \fh$.  Then a
small enough $G_{\mu}$-invariant neighbourhood $B$ of $\mu$ in the
affine plane $\mu + \fh^\circ$ is transverse to the momentum map.
($\fh^\circ$ is the annihilator of $\fh $ in $\fg^*$).  Hence
$R := \Phi \inv (B)$ is a $G_{\mu}$-invariant submanifold of $M$
containing the given \re\ $\gamma$.

We claim that in some neighbourhood of $\gamma$, $R$ is a
\emph{symplectic} submanifold of $M$.  Should the momentum map be
equivariant already with respect to the standard coadjoint action, this is a
consequence of the symplectic cross-section theorem of Guillemin and
Sternberg (cf.\ \cite{GLS}, Corollary~2.3.6).  In general, we resort
to the Witt-Moncrief decomposition described above.  As $R$ is
complementary to $T_1$ by construction, the restriction of $\omega$ to
$R$ is obtained by eliminating the second row and the second column of
$[\omega]_x$ in (\ref{eq:matrix}).  The resulting matrix is
nondegenerate, hence $R$ is symplectic in a neighbourhood of $\gamma$,
as claimed.

The action of $G_{\mu}$ on $R$ is Hamiltonian, and its momentum map is
the restriction of $\Phi$ to $R$ followed by the natural projection
$\fg^* \to \fg_{\mu}^*$.  Since $\fg^* \to \fg_{\mu}^*$ restricted to
$\mu + \fh^\circ$ is an isomorphism, the restriction $\Phi\restr{R}$
is a momentum map for the action of $G_\mu$ up to this isomorphism.
It follows that
$$
\ker d\Phi(y) = \ker d(\Phi\restr{R})(y) \qquad \forall y\in R.
$$
Moreover, because $h$ is $G$-invariant, the flow of $X_h$ preserves
the fibres of the momentum map, and so the flow preserves $R$.  It
follows that
$$
(X_h)\restr{R} = X_{(h\smallrestr{R})}.
$$

Another question that requires attending to is whether $\Phi$ inherits
openness from $\Phi\restr{R}$.  The answer is affermative in view of

\begin{Lemma}\label{lemma:slice}
  Let $K$ be a closed subgroup of a Lie group $G$ and $H$ be a closed
  subgroup of $K$.  Let $A$ be an $H$-space and $B$ a $K$-space, and
  let $f:A\to B$ be an $H$-equivariant map.  Then the map
  \begin{align*}
    F : G\times_H A &\longrightarrow G\times_K B \\
    ([g,a]_H) &\mapsto [g,f(a)]_K
  \end{align*}
  is well-defined and $G$-equivariant.  Furthermore, if $f$ is
  $H$-open, then $F$ is $G$-open.
\end{Lemma}

\begin{proof}
  The only nontrivial conclusion is the $G$-openness of $F$.  The
  diagram below commutes:
  $$
  \begin{matrix}
    G\times A & \stackrel{\mathrm{id}\times f}{\longrightarrow} & G
    \times B \cr \downarrow\hbox to0pt{$\pi_1$} & & \downarrow\hbox
    to0pt{$\pi_2$} \cr G\times_H A& \stackrel{F}{\longrightarrow} & G
    \times_K B.
  \end{matrix}
  $$
  Here $\pi_1$ and $\pi_2$ are the orbit (quotient) maps, open by
  definition of the topology on orbit spaces.  Let $U\subset G\times_H
  A$ be $G$-invariant and open.  Then $\pi_1\inv(U) = G\times U'$,
  with $U'$ open and $H$-invariant in $A$. Then $F(U)$ is open, since
  $F(U) = \pi_2\circ(\mathrm{id}\times f)(\pi_1^{-1}(U))=\pi_2(G\times
  f(U'))$.
\end{proof}

Thus, we have found a Hamiltonian sub-system $(R, G_{\mu},
\Phi\restr{R}, h\restr{R})$ for which the symmetry group $G_{\mu}$ is
compact.  Passing to this sub-system, we may and shall assume without
loss of generality that $G$ is compact and $G = G_{\mu}$.

\section{Openness of momentum maps}\label{sec:openness}

In this Section we establish Theorem~\ref{thm:openness}.  By the
result of Section~\ref{sec:reducetocompact}, we may focus our
attention on $\Phi:M\to\fg^{\ast}$, a momentum map for a Hamiltonian
action of a \emph{compact} Lie group.  Sjamaar \cite{Sj} proved that
if $\Phi$ is proper then it is $G$-open relative to its image: that
is, if $U\subset M$ is a $G$-invariant open subset, then $\Phi(U)$ is
open in $\Phi(M)$, where $\Phi(M)$ is given the subspace topology
induced from $\fg^{\ast}$.  In this section we deduce from Sjamaar's
theorem that $\Phi$ is \emph{locally $G$-open} even when it is not
proper.

Besides Lemma~\ref{lemma:slice} of Section~\ref{sec:reducetocompact},
we need two more lemmas.

\begin{Lemma} \label{lemma:polytope}
  Let $\Delta_1$ be a compact convex polytope in $\R^N$, and for $r>0$
  let
  $$\Delta_{<r} = \{a\delta\mid 0\leqslant a<r,\; \delta\in\Delta_1\}.$$
  Then $\Delta_{<r}$ is an open subset of $\Delta=\Delta_{<\infty}$,
  where the latter has the topology induced from $\R^N$.
\end{Lemma}

\begin{proof}
  There are two cases to examine depending on whether or not
  $0\in\Delta_1$. Only one ($0\in\Delta_1$) is needed for the proof of
  Theorem~\ref{thm:openness}, but we include the other for
  completeness.
  
  Case $0\in\Delta_1$.  Let $\Delta_1^j$ be the open faces of the
  polytope $\Delta_1$, and let
  $$d = \min_j\{\textrm{dist}(0,\Delta_1^j)\mid 0 \not\in
  \overline{\Delta_1^j}\}.$$
  Note that $d>0$. Then $S_d\cap\Delta_1 =
  S_d\cap\Delta$, where $S_d$ is the sphere in $\R^N$ of radius $d$.
  
  Let $\{x_n\}$ be any sequence in $\Delta$ converging to 0, and
  suppose that all $x_n\neq0$.  Write $y_n = dx_n/\|x_n\|$.  Then
  $y_n\in S_d\cap\Delta$ and so $y_n\in\Delta_1$.  Therefore $x_n =
  (\|x_n\|/d)\,y_n \in \Delta_{<r}$ provided $n$ is sufficiently large
  for $\|x_n\| < rd$ to hold.
  
  If $x_n\to x\neq 0$, we can write $x_n = a_n\delta_n$ and
  $x=a\delta$, with $\delta,\delta_n\in S_d\cap\Delta_1$.  Then
  $a_n\to a$ and $\delta_n\to\delta$.  Rescaling $\delta$ and
  $\delta_n$ if necessary the result ensues.
  
  Case $0\not\in\Delta_1$.  Let $x\in\Delta_{<r}$, so that $x =
  a\delta$ with $0\leqslant a < r$ and $\delta\in\Delta_1$.  Let
  $\{x_n\}$ be a sequence in $\Delta$ converging to $x$.  Write $x_n =
  a_n\delta_n$, with $\delta_n\in\Delta_1$. If $x=0$, then $a_n\to 0$
  as $\delta_n$ is bounded away from 0, so that $x_n\in\Delta_{<r}$
  for $n$ sufficiently large.  If on the other hand $x\neq 0$, then
  both sequences $\{a_n\}$ and $\{\delta_n\}$ are bounded and bounded
  away from 0 for sufficiently large $n$.  By rescaling $\delta_n$ if
  necessary we can arrange for $a_n$ to converge to $a$, and so
  $x_n\in\Delta_{<r}$ for $n$ sufficiently large.
\end{proof}

\begin{Lemma}  \label{lemma:quadratic}
  Let $V$ be a symplectic representation of a compact Lie group $G$,
  and let $\Phi : V\to\fh^{\ast}$ be the homogeneous quadratic momentum
  map.  Then $\Phi$ is $G$-open relative to its image.
\end{Lemma}

\begin{proof}    
  Take a $G$-invariant Hermitian metric whose imaginary part is the
  given symplectic structure on $V$.  $G$ acts as unitary
  transformations on $V$ seen as a complex vector space.
  
  Consider the unit sphere $S$ in $V$ with respect to the real part of
  the Hermitian metric, and the symplectic action of the circle group
  $\text{U}(1)$.  As $\text{U(1)}$ is the centre of the unitary group,
  the actions of $G$ and $\text{U(1)}$ commute.  The $G$-action
  descends to the symplectic manifold $\mathbb{P}V := S/\text{U(1)}$ 
  (diffeomorphic to $\mathbb{C}P^{{\rm dim}V/2-1}$).  Denote its
  momentum map by $\Phi_1 : \mathbb{P}V\to\fg^{\ast}$.  Since the
  actions of $G$ and $\text{U}(1)$ commute, $\Phi_1$ can be chosen so
  that $\Phi_1(\text{U(1)}\cdot x) = \Phi(x)$.
  
  Now let $U$ be a $G$-invariant open subset of $V$.  We want to show
  that $\Phi(U)$ is open in $\Phi(V)$, and to do so we examine two
  cases separately: (i) $0\not\in U$ and (ii) $U=B(0,\eps)$; indeed a
  general open set containing the origin is the union of sets of these
  types.
  
  Case (i): We exploit a basis for the topology of $V\setminus\{0\}
  \simeq S\times\R^+$ (diffeomorphic) consisting of `product sets'.
  Thus let $U = U_1\times(a,b)\subset S\times\R^+$, where $U_1$ is a
  $G$-invariant open subset of $S$.  Then by the homogeneity of
  $\Phi$,
  $$
  \Phi(U) = \{r^2 \nu \mid r\in(a,b),\: \nu\in \Phi(U_1)\}.
  $$
  By Sjamaar's theorem \cite{Sj, LMTW}, $\Phi(U_1)$ is open in
  $\Phi_1(\mathbb{P}V)=\Phi(S)$, and it follows that $\Phi(U)$ is open
  in $\Phi(V)$.  To see this, let $\mu\in\Phi(U)$, and let $\{\mu_n\}$
  be a sequence in $\Phi(V\setminus\{0\})$ converging to $\mu$; we may
  suppose $\mu_n\neq0$ (otherwise it is trivial).  Then $\mu = \Phi(v,
  r) = r^2\Phi_1(v)$ for some $(v, r)\in S\times (a,b)$.  Since
  $\mu/r^2 \in \Phi_1(\mathbb{P}V)$, there is a sequence $r_n\to r$
  satisfying, $\mu_n/r_n^2 \in \Phi_1(\mathbb{P}V)$ for all $n$.
  $\Phi_1$ being open by Sjamaar's theorem, there is a sequence
  $(v_n)$ in $S$ converging to $v$ such that $\mu_n/r_n^2 =
  \Phi_1(v_n)$; in other words $\mu_n = \Phi(v_n, r_n)$, and $(v_n,
  r_n) \to (v, r)$.  Consequently $\Phi(U)$ is open, as required.

     Case (ii):  Let $U=B(0,\eps)$, the open ball in $V$ with centre
  $0$ and radius $\eps$.  Because $G$ is compact, we can and do identify 
  $\fg$ with $\fg^{\ast}$, and the adjoint action with the coadjoint action.  
  Let $\ft^+$ be a positive Weyl chamber in $\fg=\fg^{\ast}$, and let
  \begin{align*}
  \Delta_1 &= \Phi_1(\mathbb{P}V) \cap \ft^+, \\
  \Delta   &= \Phi(V) \cap \ft^+.  
  \end{align*}
  By the convexity theorem of Atiyah-Guillemin-Sternberg-Kirwan,
  $\Delta_1$ is a convex polytope, and by the homogeneity of $\Psi$,
  $\Delta=\R^+\Delta_1$.  $U$ is $G$-invariant, and
  $$
  \Phi(U)\cap \ft^+ = \Delta_{<\eps^2} = \bigcup_{0\leqslant
    r<\eps^2}r\Delta_1.
  $$
  By Lemma \ref{lemma:polytope}, $\Delta_{<\eps^2}$ is open in $\Delta$.

     To finish the proof that $\Phi(U)$ is open in $\Phi(V)$, note 
  that $U$ is $G$-invariant, so that both $\Phi(U)$ and $\Phi(V)$ are 
  $G$-invariant subsets of $\fg^{\ast}$. Since $\ft^+$ and $\fg^{\ast}/G$ 
  are homeomorphic, we have that $\Phi(U)/G$ is open in $\Phi(V)/G$, and 
  the result ensues.
\end{proof}

\noindent
\emph{Remark} This lemma illustrates why it is important to consider
$G$-openness rather than openness, for momentum maps are not in
general open.  For example the momentum map for the ${\rm SO}(3)$
action on $T^*\R^3$, namely $\Phi(q,p) = q\times p$, is not open. For
example, the images of sufficiently small neighbourhoods of $(q,p) =
(e_1,e_1)$ are not neighbourhoods of $0$: in particular they do not
contain nonzero points of the $e_1$-axis.

\medskip

Equipped with these lemmas, we are in a position to prove
Theorem~\ref{thm:openness}.

\begin{proof} 
  Let $x\in M$ and $\mu=\Phi(x)$, and let $U_1$ be the $G$-invariant
  neighbourhood of $x$ whose existence is guaranteed by the
  Marle-Guillemin-Sternberg normal form
  (Section~\ref{sec:reducetocompact}).  Let $V_1$ be a $G$-invariant
  tubular neighbourhood of $\mu$, so that
  $$
  V_1 \simeq G \times_{G_\mu} \cO,
  $$
  where $\cO$ is a neighbourhood of 0 in the slice $\fg_\mu^\sharp$
  (the fixed-point set of the action of the centre of $G_{\mu}$ on
  $\fg^{\ast}$, which is isomorphic to $\fg_\mu^*$).  Finally, let $U
  = U_1\cap\Phi^{-1}(V_1)$. Thus, as symplectic $G$-spaces, $U\simeq
  G\times_{G_x} (\fm^*\times Y)$, where $Y$ is a $G_x$-invariant
  neighbourhood of $0$ in $N_1$, and it therefore suffices to show
  that the momentum map (\ref{eq:MGS}) is open.  Since $G$ is compact,
  we can take the cocycle $\theta$ in (\ref{eq:MGS}) to vanish.
  
  By Lemma~\ref{lemma:quadratic}, the quadratic momentum map
  $\Phi_{G_x}:N_1\to\fg_x^*$ is open relative to its image.  It
  follows that the restriction of $\Phi$ to the slice
  $\fm^{\ast}\times Y \to \fg_\mu^*$ is open relative to \emph{its}
  image.  Lemma~\ref{lemma:slice} now applies, with $f$ replacing the
  restriction of $\Phi$ to the slice and $F$ replacing $\Phi$.
\end{proof}


\section{Persistence of extremal relative equilibria}\label{sec:persistence}

In this section we establish Theorem~\ref{thm:persistence}, using
Theorem~\ref{thm:openness} which was proved in
Section~\ref{sec:openness}.  By the results of
Section~\ref{sec:reducetocompact} we may assume $G$ to be compact.

We treat the case when $\gamma$ is a minimal \re; the maximal case is
similar. Let $U_0$ be the $G$-invariant neighbourhood of $\gamma$
guaranteed by Theorem~\ref{thm:openness}.  The minimality of $\gamma$
means that there is a precompact $G$-invariant neighbourhood $U\subset U_0$ 
of $\gamma$ such that
$$
h\restr{\Phi\inv(\mu) \cap \overline{U}} \geqslant h(\gamma) \quad
\text{with equality only on} \quad G\cdot \gamma \cap U.
$$

Suppose Theorem~\ref{thm:persistence} is false.  Let $\{\mu_n\}$ be a
sequence of points in $\Phi(U)$ (which is open in $\Phi(U_0)$)
converging to $\mu$, such that the restriction of $h$ to
$\Phi^{-1}(\mu_n)\cap U$ has no minimum.  However, by compactness, the
restriction of $h$ to $\overline{\Phi^{-1}(\mu_n)\cap U}$ has a
minimum, say at $y_n\in\overline{U}\diagdown U$.  Also by compactness,
$y_n\to y$, with $y\in \overline{U}\diagdown U$ (possibly after
passing to a subsequence).

We claim that there is a sequence $\{x_n\}$ converging to some $x\in
G\cdot\gamma$, with $\Phi(x_n)=\mu_n$.  Granted that claim,
we have $h(x) =
h(G\cdot\gamma) < h(y)$ by construction of $y$.  On the other hand,
for each $n$, $h(y_n) < h(x_n)$.  In the limit we get $h(y)\leqslant h(x)$,
which is a contradiction.

The existence of the sequence $\{x_n\}$ is a consequence of the
openness property.  Indeed, we can choose a nested sequence $\{U_n\}$ of
$G$-invariant neighbourhoods of $G\cdot\gamma$ whose intersection is
$G\cdot\gamma$, such that $\mu_n\in\Phi(U_n)$.  Choosing $x_n\in
\Phi^{-1}(\mu_n)\cap U_n$ gives (after passing to a subsequence if
necessary) a sequence converging to a point $x\in G\cdot\gamma$, as
claimed.


\section{An example: the affine rigid body}\label{sec:ARB}

The problem of affine rigid bodies (\emph{alias} Riemann ellipsoids)
has a long and important history, dating back perhaps to when Newton
correctly suggested that the Earth was an oblate spheroid.  Since
then, it has been studied by such illustrious figures as Maclaurin,
Jacobi, Dirichlet, Riemann, and Poincar\'e.  A classical discussion can
be found in the book of Chandrasekhar \cite{Chand}; for a recent
account from the symmetry perspective, we refer to Roberts and Sousa
Dias \cite{RSD99}.

An affine rigid body models a mass of ideal fluid evolving in time in
such a manner that it always remains an ellipsoid.  This is a
Hamiltonian system whose configuration space is either $\cQ
= \text{SL}(3, \R)$ or $\text{GL}(3, \R)$ depending on whether one is
modelling incompressible or compressible fluids.  The matrix $Q \in
\cQ$ represents the configuration that is the image of a sphere under
$Q$, an ellipsoid whose semi-axes are given by the
singular values of $Q$.  It is supposed that the potential energy
depends only on the shape of the ellipsoid, and so is invariant under
the symmetry group $G=\text{SO}(3)\times\text{SO}(3)$, the first copy
of $\text{SO}(3)$ acting by multiplication on the left, and the second
by multiplication on the right.

The phase space is then the cotangent bundle $T^{\ast}\cQ$, and the
group $G$ acts by cotangent lift on $T^{\ast}\cQ$.  Accordingly, the
momentum map has two `components'
$$
\Phi_L,\Phi_R : T^{\ast}\cQ \to \so(3)^{\ast}
$$
given by
$$
\Phi_L(Q,P)=\frac{1}{2}\left(PQ^T-QP^T\right) ,\qquad
\Phi_R(Q,P)=\frac{1}{2}\left(P^TQ-Q^TP\right) .
$$

The particular example of the potential energy function used by
Dirichlet, Riemann and others is the \emph{self-gravitating}
potential.  Other potentials arise in linear elasticity theory.  In
most of these examples the potential energy $V(Q)$ has a minimum at
the round sphere $Q=I$.  It then follows that this point is an
equilibrium and indeed an \emph{extremal} (relative) equilibrium.
From Theorem~\ref{thm:persistence} we deduce

\begin{Corollary}
  Suppose the potential energy has a minimum at the point $(I,0) \in
  T^{\ast}\cQ$. Then there exist $\eps_L, \eps_R>0$ and a
  $G$-invariant neighbourhood $U$ of $(I,0)$ in $T^{\ast}\cQ$ such
  that for all $(\mu_L, \mu_R) \in \so(3)^*\times\so(3)^*$ with
  $\|\mu_L\| < \eps_L$, $\|\mu_R\|<\eps_R$ there is an extremal
  relative equilibrium of the affine rigid body in $U$ with momentum
  $(\Phi_L, \Phi_R)=(\mu_L, \mu_R)$.
\end{Corollary}

Recall \cite{Mo97} that extremal relative equilibria are
Lyapunov-stable relative to $G$, and by \cite{LS} they are then
Lyapunov-stable relative to $G_\mu$ as well.


\section{A counter-example: plane point vortices}\label{sec:vortices}

In this section we give an example to the effect that the hypothesis
on the compactness of the momentum isotropy subgroup $G_\mu$ is
necessary in Theorem~\ref{thm:persistence}.  Consider the symplectic
manifold
$$
M = \C^{\,N} \, \diagdown \bigcup_{k \ne l} \{z_k = z_l\}, \quad
\omega = \frac{i}{2} \sum_{k=1}^N \Gamma_k dz_k\wedge d\overline{z}_k,
\quad \Gamma_1, ~\ldots, \Gamma_N \in \R 
$$
on which the Euclidean group $G = {\rm SE}(2) = \R^2 \rtimes
\text{SO}(2)$ acts diagonally.  This action is free, proper, and
Hamiltonian, and has a momentum map $\Phi : M \to \fg^{\ast}$.  We use
the identification $\fg^{\ast} \simeq \C \times \R$ and denote the
components of $\Phi$ by
$$
(\prc, \prr) : (z_1, ~\ldots, z_N) \mapsto 
\left( i\sum_{k=1}^N \Gamma_k z_k,
\sum_{k=1}^N \Gamma_k \frac{|z_k|^2}{2} \right).
$$
It can be shown that $\Phi$ is equivariant with respect to the
standard (unmodified) coadjoint action if and only if 
$\sum_{k=1}^N \Gamma_k = 0$.  As the $G$-invariant
Hamiltonian we take
$$
h(z_1, ~\ldots, z_N) = -\frac{1}{2\pi}\sum_{k<l}
\Gamma_k\Gamma_l \log |z_k - z_l|.$$
Hamilton's equation reads
$$
\frac{dz_k}{dt} = 
\frac{2}{i}\frac{\partial h}{\partial (\Gamma_k\overline{z}_k)} =
-\frac{1}{2\pi i}\sum_{l\neq k}
     \frac{\Gamma_l}{\overline{z}_k - \overline{z}_l}
\qquad (k = 1, ~\ldots, N).
$$
This system describes the motion of $N$ interacting plane point 
vortices with vorticities $\Gamma_1, ~\ldots, \Gamma_N$.
See for example \cite{Aref83}.

We study the case of 3 vortices with vorticities 1, 1, $-2$.  Let us
call the \emph{axis} the subset $0 \times \R^{+}$ of $\C\times\R\simeq
\fg^{\ast}$.  A theorem of Synge, \cite[Theorem 6]{Synge49}, tells us
that the relative equilibria $\gamma = (z_1, z_2, z_3)$ for $h$ are of
two types:
\begin{enumerate}
\item $\Phi(\gamma)$ is on the axis, i.e.\ $\prc(\gamma) = 0$, in
  which case $z_1, z_2, z_3$ are collinear, with $z_3$ midway between
  $z_1$ and $z_2$;
\item $\Phi(\gamma)$ is off the axis, i.e.\ $\prc(\gamma) \ne 0$, in
  which case $z_1, z_2, z_3$ form an equilateral triangle.
\end{enumerate}

\noindent In an equilateral relative equilibrium (type 2), $z_2 - z_3 =
e^{\pm i\pi/3}(z_1 - z_3)$, from which we calculate easily that
$$
e^{2\pi h(\gamma)} = \frac{1}{3\sqrt{3}}|\prc(\gamma)|^3.
$$

\begin{Proposition}\label{prop:vortices} 
  In a system of 3 vortices with vorticities 1, 1, $-2$, let $\gamma$
  be a collinear \re\ for $h$ (type 1), $U$ the $G$-invariant
  neighbourhood of $\gamma$ in $M$ defined by $e^{2\pi h(U)} > e^{2\pi
    h(\gamma)}/3$, and $V = (D\times\R\diagdown$ axis) a punctured
  neighbourhood of $\mu = \Phi(\gamma)$ in $\fg^{\ast}$ where
  $D\subset \C$ is the disc of radius $(\sqrt{3}\, e^{2\pi
    h(\gamma)})^{1/3}$ centred at $0$.  Then $\gamma$ is extremal, but
  $\Phi\inv(V)\cap U$ contains no \re\ for $h$.
\end{Proposition}

\begin{proof}
  Since the action of $G$ is free, $\Phi$ is a submersion, and
  $\operatorname{codim}(\Phi\inv(\mu)) =\operatorname{codim}(\mu) = 3
  = \dim(G) = \dim(G_{\mu})$.  This means that
  $\Phi^{-1}(\mu)/G_{\mu}$ is discrete (by explicit calculation, in
  fact a single point), hence $\gamma$ is trivially extremal.  Suppose
  a \re\ $\gamma'$ exists in $\Phi\inv(V)\cap U$.  Since
  $\Phi(\gamma')$ is off the axis, $\gamma'$ is an equilateral
  triangle (type 2) and $e^{h(\gamma')} =
  \frac{1}{3\sqrt{3}}|\prc(\gamma')|^3 <
  \frac{1}{3\sqrt{3}}|\sqrt{3}\, e^{h(\gamma)}|$.  This is
  incompatible with $e^{h(\gamma')} > e^{h(\gamma)}/3$.
\end{proof}

\begin{Remark}
  In Proposition~\ref{prop:vortices}, $M, G, \Phi, h, \gamma$ satisfy all
  the hypotheses in Theorem~\ref{thm:persistence} except the
  compactness of the isotropy group of $\Phi(\gamma)$, which is $G =
  {\rm SE}(2)$ itself.  The failure of $\gamma$ to persist to nearby
  levels of the momentum map shows that this compactness hypothesis is
  essential.  As observed in the proof, the momentum map is a
  submersion and the conclusion of Theorem~\ref{thm:openness} is still
  true.
  
  Note, however, that $\gamma$ persists when $\mu$ is perturbed along
  the axis $0\times \R$, which is the annihilator of the noncompact
  part $\R^2$ of the Lie algebra and so in a natural way the dual of
  the Lie algebra of the compact subgroup $\text{SO}(2)$ of $G$.  This
  partial persistence in `compact directions' is in fact also covered
  by Theorem~\ref{thm:persistence}, via reduction by stages and
  Corollary~\ref{cor:rbys} below.
  
  It is interesting that, in the model of $N$ plane vortices, if the
  sum of the vorticities does not vanish, then the momentum map is
  equivariant with respect to a modified coadjoint action (with a
  nontrivial cocycle $\theta$), and for this modified coadjoint action
  all the isotropy subgroups are isomorphic to $\text{SO}(2)$ and so
  are compact.  Hence Theorem~\ref{thm:persistence} implies
  that the \re\ persists to nearby values of momentum.  
  This does not contradict the theorem of Synge quoted above, for 
  in the case of nonvanishing total vorticity all nearby values of 
  momentum are realisable by collinear configurations of 3 vortices. 

\end{Remark}


\section{Reduction by stages}\label{sec:stages}

Let $K$ be a normal subgroup of a Lie group $G$ with quotient $L=G/K$.
Roughly speaking, we say that reduction by stages works if reduction
by $G$ coincides with reduction first by $K$ and then by $L$.  In
detail, reduction by stages describes the following general procedure.

At the level of Lie algebras and their duals, we have $\fk\subset\fg$
and $\fl^*\subset\fg^*$.  Moreover, the inclusion $\fl^*
\hookrightarrow \fg^*$ naturally identifies $\fl^*$ with the
annihilator $\fk^\circ$ of $\fk$ in $\fg^*$.  Let $\pi:\fg^*\to\fk^*$
be the canonical projection.  

Suppose $G$ acts in a Hamiltonian manner on a symplectic manifold $M$
with momentum map $\Phi_G:M\to\fg^*$.  This action restricts to an
action of $K$, and the momentum map is just $\Phi_K=\pi\circ\Phi_G :
M\to\fk^*$.  $G$ and $K$ act on $\fg^*$ and $\fk^*$ in such a way as
to make the momentum maps equivariant.  For $\mu\in\fg^*$, we shall
use the notation $M\sympquot{\mu}{G}$ to denote the reduced space
$\Phi_G^{-1}(\mu)/G_\mu$.

The modified coadjoint action $\Coad^\theta$ of $G$ on $\fg^*$
descends to an action on $\fk^*$, and so $L$ acts in a natural way on
the set of $K$-orbits in $\fk^*$.  Let $L_\nu$ be the subgroup of $L$
that preserves the coadjoint orbit $K\cdot\nu$; one can show that
$L_\nu\simeq G_\nu/K_\nu$.  It follows that $L_\nu$ acts on
$M\sympquot{\nu}{K}$ in a natural way, preserving the symplectic
structure (compare \cite{MR00}).

The affine subspace $\pi\inv(\nu)$ of $\fg^*$ can be identified with
$\fk^\circ$, and hence with $\fl^*$, by translation in $\fg^*$.  Let
$\rho_\nu:\pi\inv(\nu)\to \fl^*_\nu$ be the composite of such an
identification with the natural projection $\fl^*\to\fl_\nu^*$. Define
$$\Phi_\nu : M\sympquot{\nu}{K} \to \fl_\nu^*$$%
by first restricting $\Phi_G$ to $\Phi_K\inv(\nu)$ (whose values lie
in $\pi\inv(\nu)$), passing to the quotient $M\sympquot{\nu}{K}$ and
finally applying the identification $\rho_\nu$.  If this map
$\Phi_\nu$ is well-defined and is a momentum map for the action of
$L_\nu$, then we say \emph{reduction by stages works\/} in this
context provided in addition that
$$M\sympquot{\mu}{G} \simeq \left(M\sympquot{\nu}{K}\right)
\sympquot{\sigma}{L_\nu}$$ %
at least at the level of connected components, where $\nu=\pi(\mu)$
and $\sigma=\rho_\nu(\mu) = \Phi_\nu(\Phi_G^{-1}(\mu))$.  The
isomorphism between the two spaces should be as needed in the context;
here a homeomorphism is sufficient, though more generally one might
require an isomorphism as symplectic stratified spaces \cite{SL91}.
The papers \cite{MMPR98,MMOPR00} explain the current state of the art
on reduction by stages for free actions.

If the action of $K$ is free, then the partially reduced space
$M\sympquot{\nu}{K}$ is a smooth symplectic manifold, and one can
apply Theorem~\ref{thm:persistence} to the resulting $L_\nu$-invariant
system, provided $(L_\nu)_\sigma$ is compact.  The result is that if
the \re\ $\gamma$ in question is extremal in $M\sympquot{\nu}{K}$ then
it persists to nearby values of the $L_\nu$-momentum map $\Phi_\nu$.

It often happens that $L_\nu=L$, in which case $\rho_\nu$ is just a
translation of $\mu+\fk^\circ$ to $\fk^\circ$, and one can ask for
persistence to all $\mu'\in V\subset (\mu+\fk^\circ)$.


\begin{Corollary} \label{cor:rbys}
  Let $G,M,\Phi,K,L,\mu,\nu$ and $\sigma$ be as above, with $K$ acting
  freely on $M$, and such that reduction by stages works.  Suppose
  that $L_\nu=L$, and that $L_\sigma$ is compact.  Let $h$ be a
  $G$-invariant Hamiltonian on $M$ for which $\gamma\subset M$ is an
  extremal relative equilibrium with $\mu=\Phi(\gamma)$.  Let $U_0$ be
  the $L$-invariant neighbourhood of $\gamma$ in $M\sympquot{\nu}{K}$
  guaranteed by Theorem~\ref{thm:openness}. Then there is an
  $L$-invariant neighbourhood $V$ of\/ $\mu$ in
  $\Phi(U_0)\cap\fk^\circ$ such that for every $\sigma'\in V$ there is
  a relative equilibrium in $\Phi^{-1}(\sigma')\cap U_0$.
\end{Corollary}

For the vortex system of Section~\ref{sec:vortices} the corollary
applies as follows: $G={\rm SE}(2)$, $K=\R^2$ the normal subgroup of
${\rm SE}(2)$ consisting of translations, $\mu=(0,r)$, and $\nu=0$.
Then $L_\nu = L = {\rm SE}(2)/\R^2\simeq{\rm SO}(2)$, which is
compact.  The subspace $\fk^\circ$ is a line (the `axis' of
Section~\ref{sec:vortices}) consisting of the coadjoint orbits that
are isolated points.

\begin{Remark}
  The freeness hypothesis on the $K$-action could easily be relaxed to
  local freeness.  In the general setting where the action is not
  locally free, the same argument can be applied to the symplectic
  stratum in $M\sympquot{\nu}{K}$ containing the image of $\gamma$.
\end{Remark}


\bigskip

\noindent\textbf{Acknowledgement}.
We would like to thank Eugene Lerman for many useful discussions, and
for pointing out an error in a preliminary draft.


 \end{document}